\documentclass[smallextended,envcountsect,]{svjour3}
\smartqed
\usepackage{graphicx}
\usepackage{amsmath, amsfonts, amssymb, mathtools}
\usepackage{authblk, color, bm, graphicx, epstopdf, url}
\usepackage[font = small, labelfont = bf, labelsep = period]{caption}
\usepackage{subcaption}
\usepackage{xspace}
\usepackage[dvipsnames]{xcolor}
\usepackage[framemethod=Tikz]{mdframed}
\usepackage{enumitem,hyperref}
\usepackage{booktabs, tabularx, multirow}
\usepackage{algorithm2e}

\allowdisplaybreaks

\date{\today}

\begin{document}

\title{Condensed interior-point methods: porting reduced-space approaches on GPU hardware}

\authorrunning{Pacaud et al.}
\author{
  François Pacaud \and Sungho Shin \and Michel Schanen \and Daniel Adrian Maldonado \and Mihai Anitescu
  \thanks{Mihai Anitescu dedicates this work to the 70-th birthday of Florian Potra. Florian, thank you for the great contributions to optimization in general, and interior point methods in particular, and for initiating me and many others in them.}
}

\institute{All the authors are affiliated to Argonne National Laboratory.
}
\maketitle

\abstract{
  The interior-point method (IPM) has become the workhorse
  method for nonlinear programming. The performance of IPM is
  directly related to the linear solver employed to factorize the Karush--Kuhn--Tucker (KKT) system
  at each iteration of the algorithm.
  When solving large-scale nonlinear problems, state-of-the art IPM solvers
  rely on efficient sparse linear solvers to solve the KKT system.
  Instead, we propose a novel reduced-space
  IPM algorithm that condenses the KKT system into a dense matrix
  whose size is proportional to the number of degrees of freedom in the problem.
  Depending on where the reduction occurs we derive two variants of the reduced-space
  method: linearize-then-reduce and reduce-then-linearize. We adapt their
  workflow so that the vast majority of computations are accelerated on GPUs.
  We provide extensive numerical results
  on the optimal power flow problem, comparing our GPU-accelerated reduced
  space IPM with Knitro and a hybrid full space IPM algorithm.
  By evaluating the derivatives on the GPU and solving the KKT system
  on the CPU, the hybrid solution is already significantly faster than the CPU-only solutions.
  The two reduced-space algorithms go one step further by solving the KKT system entirely
  on the GPU.
  As expected, the performance of the two reduction algorithms depends
  intrinsically on the number of available degrees of freedom:
  their performance is poor when the problem has many degrees of freedom,
  but the two algorithms are up to 3 times faster than Knitro as soon as the
  relative number of degrees of freedom becomes smaller.


\section{Introduction}

Most optimization problems in engineering consist on minimizing a
cost subject to a set of equality constraints that represent the
physics of the problem. Often, only a subset of the problem variables are
actionable. A particular instance of these types of problems is the
Optimal Power Flow (OPF), which formulates as a large-scale nonlinear nonconvex
optimization problem. It is one of the critical power system analysis
carried out multiple times per day by any electrical power system market
operator~\cite{frank2012optimal}. In brief, it selects the optimal real power of
dispatchable resources (typically power plants) subject to
physical constraints (the power flow constraints, or power balance), and operational constraints
(e.g voltage or line flow limits)~\cite{cain2012history}.
It is  challenging to solve to optimality, particularly since its solution is
needed within a prescribed and fairly tight time limit. Since the 1990s, the efficient solution
of OPF has relied on the interior-point method (IPM).
In particular, state-of-the-art numerical tools developed for sparse IPMs can efficiently handle the sparse structure of typical OPF problems.
At each iteration of the IPM algorithm, the descent direction
is computed by solving a Karush--Kuhn--Tucker (KKT) system, which requires the factorization
of large-scale ill-conditioned unstructured symmetric indefinite matrices \cite{nocedal_numerical_2006}.
For that reason, current state-of-the-art OPF solvers combine
a mature IPM solver~\cite{wachter2006implementation,waltz2006interior}
together with a sparse Bunch--Kaufman factorization routine~\cite{duff2004ma57,schenk2004solving}.
Along with an efficient evaluation of the derivatives, this method
is able to solve efficiently OPF problems with up to 200,000 buses on
modern CPU architectures~\cite{kardos2018complete}.

\subsection{Interior-point method on GPU accelerators: state of the art}
Most upcoming HPC architectures are GPU-centric, and we can leverage
these new parallel architectures to solve very-large OPF instances.
However, porting IPM to the GPU is nontrivial because GPUs are based
on a different programming paradigm from that of CPUs: instead of computing
a sequence of instructions on a single input (potentially dispatched on different threads or processes),
GPUs run the same instruction simultaneously on hundreds of threads (SIMD paradigm:
\emph{Single Instruction, Multiple Data}).
Hence, GPUs shine when the algorithm can be decomposed into simple instructions
 running entirely in parallel where the same instruction is executed on different data in lockstep (as is the case for most dense linear algebra).
In general, not all algorithms are fully amenable to this paradigm. One notorious example is branching in the control flow:
when the instructions have multiple conditions,
dispatching the operations on multiple threads makes execution in lockstep impossible.

Unfortunately, factorization of unstructured sparse indefinite matrix is  one of these edge cases:
unlike dense matrices, sparse matrices have unstructured sparsity, rendering most sparse algorithms difficult to parallelize.
Thus, implementing a sparse direct solver on the GPU is nontrivial,
and the performance of current GPU-based sparse linear solvers lags far behind that of their CPU equivalents~\cite{swirydowicz2021linear,tasseff2019exploring}.
Previous attempts to solve nonlinear problems on the GPU
have circumvented this problem by relying on
iterative solvers \cite{cao2016augmented,schubiger2020gpu}
or on decomposition methods~\cite{kim2021leveraging}.
Here, we have chosen instead to revisit the original reduced-space algorithm
proposed in~\cite{dommel_optimal_1968}: this method
\emph{condenses} the KKT system into a dense matrix, whose
size is small enough to be factorized efficiently on the GPU with dense direct linear algebra.


\subsection{Reduced-space interior-point method}
Reduced-space algorithms have been studied for a long time.
In~\cite{dommel_optimal_1968} the authors introduced a
reduced-space algorithm, one of the first effective methods
to solve the OPF problem. The method
has several first-order variants, known as the generalized reduced-gradient
algorithm \cite{abadie1969carpentier} or the gradient
projection method~\cite{rosen1961gradient}, whose theoretical
implications are discussed in \cite{gabay1982minimizing}.
The extension of the reduced-space method to second-order comes later
\cite{sargent1974reduced}---as well as its application to OPF \cite{burchett1984quadratically}---the method becoming during the 1980s a particular case of the
sequential quadratic programming algorithm~\cite{coleman1982nonlinear,fletcher2013practical,nocedal1985projected,gurwitz1989sequential}.
However, the (dense) reduced Hessian has always been challenging
to form explicitly, favoring the development of approximation
algorithms for second-order derivatives, based on
quasi-Newton~\cite{biegler1995reduced} or on
Hessian-vector products~\cite{biros_parallel_2005}.
The reduced-space method was extended to IPM  in the late 1990s
\cite{cervantes2000reduced} and was already adopted  to solve OPF~\cite{jiang2010reduced}. We refer to \cite{kardos2020reduced} for a recent report
describing the application of reduced-space IPM to OPF.

\subsection{Contributions}

Our reduced-space algorithm is built on this  extensive previous work.
In Section~\ref{sec:reduction} we propose a tractable reduced-space IPM algorithm, allowing
 the OPF problem to  be solved entirely on the GPU.
Instead of relying on a direct sparse solver,
our reduced-space IPM \emph{condenses} the KKT system into a dense linear system whose size  depends only on the number of control variables in the reduced problem (here the number of generators).
When the number of degrees of freedom is much smaller than the total number of variables,
the KKT system size can be dramatically reduced.
In addition, we establish a formal connection
between the reduced IPM algorithm and the seminal reduced-gradient method
of Dommel and Tinney~\cite{dommel_optimal_1968}.
Depending on whether the reduction occurs at the KKT system level or directly at
the nonlinear level, we propose two different reduced-space algorithms:
\emph{Linearize-then-reduce} and \emph{Reduce-then-linearize}.
We describe a parallel implementation of the reduction algorithm
in Section~\ref{sec:streamline} and show how we can exploit
efficient linear algebra kernels on the GPU to accelerate the algorithm.
We discuss an application of the proposed method to the OPF problem
in Section~\ref{sec:results}:
our numerical results show that both the reduced IPM and its feasible
variant are able to solve large-scale OPF instances---with up to 70,000 buses---entirely on the GPU.
This result improves on the previous results reported in \cite{kardos2020reduced,pacaud2021feasible}.
As expected, the reduced-space algorithm is competitive when the problem has fewer degrees of
freedom, but it achieves respectable performance (within a factor of 3 compared with state-of-the-art methods) even
on the less favorable instances.
To the best of our knowledge, this is the first time a second-order GPU-based NLP solver
matches the performance of state-of-the art CPU-based solvers on the resolution of OPF problems.

\section{Reduced interior-point method}
\label{sec:reduction}

In \S\ref{sec:reduced:formalism} we introduce the problem formulation under consideration.
In this problem the independent variables (the \emph{control}, associated to the problem's \emph{degrees}
of freedom) are split
from the dependent variables (the \emph{state}). The reduction is akin to a Schur complement reduction
and can occur at the linear algebra or  the nonlinear levels. In \S\ref{sec:reduced:ripm}
we present our first method, \emph{linearize-then-reduce}, performing the reduction
directly on the KKT system.
In \S\ref{sec:reduced:fripm} we show that we can reduce equivalently
the nonlinear model using the implicit function theorem, giving
our second method: \emph{reduce-then-linearize}.

\subsection{Formalism}
\label{sec:reduced:formalism}
The problem we  study in this section has a particular
structure: the optimization variables are divided into
a control variables $\bm{u} \in \mathbb{R}^{n_u}$ and a state variables
$\bm{x} \in \mathbb{R}^{n_x}$. The control and the state variables
are coupled together via a set of equality constraints:
\begin{equation}
  \label{eq:stateequation}
  g(\bm x, \bm u) = 0 \, ,
\end{equation}
with $g:\mathbb{R}^{n_x} \times \mathbb{R}^{n_u} \to \mathbb{R}^{n_x}$
(note that we choose the dimension of the output space to be equal to $n_x$, the dimension of the state).
The function $g$ is often related to the physical
equations of the problem; and,
depending on the applications, it can encode balance equations (optimal power flow, the primary object of study in this work), but also discretizations of dynamics (optimal control),
or partial differential equations (PDE-constrained optimization).

We call Equation~\eqref{eq:stateequation} the \emph{state equation} of the problem.
Further, an objective function $f:\mathbb{R}^{n_x} \times \mathbb{R}^{n_u} \to \mathbb{R}$
and a set of generic constraints $h:\mathbb{R}^{n_x} \times \mathbb{R}^{n_u} \to \mathbb{R}^m$ are given.
The optimization problem in state-control form can be expressed as
\begin{equation}
  \label{eq:originalproblem}
   \min_{\bm{x},\bm{u}} \; f(\bm{x}, \bm{u}) \quad
\text{subject to} \quad
  \left\{
  \begin{aligned}
    & \bm{u} \geq 0 \; , \quad \bm{x} \geq 0 \; , \\
    & g(\bm{x}, \bm{u}) = 0 \;, \quad h(\bm{x}, \bm{u}) \leq 0  \, .
  \end{aligned}
  \right.
\end{equation}
The focus of this paper is the efficient solution of~\eqref{eq:originalproblem}
using the IPM. For that reason one often prefers to introduce
slack variables $\bm{s} \in \mathbb{R}^m$ for the nonlinear inequality constraints.
Then \eqref{eq:originalproblem} can be rewritten as
\begin{equation}
  \label{eq:slackedproblem}
    \min_{\bm{x},\bm{u}, \bm{s}} \; f(\bm{x}, \bm{u}) \quad
    \text{subject to} \quad
    \left\{
  \begin{aligned}
     \quad & \bm{u} \geq 0 \; , \quad \bm{x} \geq 0 \; , \quad \bm{s} \geq 0 \\
                            & g(\bm{x}, \bm{u}) = 0 \; , \quad h(\bm{x}, \bm{u}) + \bm{s} = 0  \, .
  \end{aligned}
  \right.
\end{equation}
The Lagrangian of the problem~\eqref{eq:slackedproblem} is
\begin{multline}
  \label{eq:lagrangian}
  L(\bm{x}, \bm{u}, \bm{s}; \bm{\lambda}, \bm{y}, \bm{v}, \bm{w}, \bm{z}) =
  f(\bm{x}, \bm{u}) + \bm{\lambda}^\top g(\bm{x}, \bm{u}) \\
  + \bm{y}^\top \big(h(\bm{x}, \bm{u}) + \bm{s}\big)
  - \bm{v}^\top \bm{x} - \bm{w}^\top \bm{u} - \bm{z}^\top \bm{s} \; ,
\end{multline}
where $\bm{\lambda} \in \mathbb{R}^{n_x}$ is the multiplier associated with the state equation~\eqref{eq:stateequation},
$\bm{y} \in \mathbb{R}^{m}$ the multiplier associated with the inequality constraints $h(\bm x, \bm u) + \bm s = 0$,
and $\bm{v} \in \mathbb{R}^{n_x}, \bm{w} \in \mathbb{R}^{n_u}, \bm{z} \in \mathbb{R}^m$
the multipliers associated with the bound constraints.

In what follows, we assume that $f$, $g$, and $h$ are twice continuously differentiable
on their domain. Throughout the article we denote
\begin{equation*}
  \begin{aligned}
    & \bm{g} = \nabla_{(\bm x, \bm u)} f(\bm x, \bm u) \in \mathbb{R}^{n_x+n_u}& \text{gradient of the objective} \\
    & A = \partial_{(\bm x, \bm u)} h(\bm x, \bm u) \in \mathbb{R}^{m \times (n_x+n_u)}& \text{Jacobian of the inequality cons.} \\
    & G = \partial_{(\bm x, \bm u)} g(\bm x, \bm u) \in \mathbb{R}^{n_x \times (n_x+n_u)}& \text{Jacobian of the equality cons.} \\
      & W = \nabla_{(\bm x, \bm u)}^2 L(\bm x, \bm u, \bm s; \bm{\lambda}, \bm{y}, \bm{v}, \bm{w}, \bm{z})
    & \text{Hessian of Lagrangian.}
  \end{aligned}
\end{equation*}
We partition the first- and second-order derivatives into the blocks associated with
the state variables $\bm x$ and the control variables $\bm u$; that is,
\begin{equation*}
  \bm{g} = \begin{bmatrix}
   \bm{g}_u \\ \bm{g}_x
  \end{bmatrix}
  \, , \quad
  \begin{aligned}
    & A = \begin{bmatrix}
      A_u & A_x
    \end{bmatrix} \\
    & G = \begin{bmatrix}
      G_u & G_x
    \end{bmatrix}
  \end{aligned}
  \, , \quad \text{and} \quad
  W = \begin{bmatrix}
    W_{uu} & W_{ux} \\
    W_{xu} & W_{xx}
  \end{bmatrix}
  \; .
\end{equation*}

\subsection{Linearize-then-reduce}
\label{sec:reduced:ripm}
Now we discuss the first reduction method: linearize-then-reduce.
This method exploits the structure
of the problem~\eqref{eq:originalproblem} directly at the linear algebra level.

\subsubsection{Successive reductions of the KKT system}
We first present the KKT system in an augmented form and show how we can
reduce it by  removing from the formulation first the equality constraints and then the inequality constraints.

\paragraph{KKT conditions.}
The KKT conditions associated with the standard formulation~\eqref{eq:slackedproblem} are
\begin{subequations}
  \label{eq:kkt:fullspace}
\begin{align}
  & \bm{g}_u + G_u^\top \bm{\lambda} + A_u^\top \bm{y} - \bm{w} = 0  , \\
  & \bm{g}_x + G_x^\top \bm{\lambda} + A_x^\top \bm{y} - \bm{v} = 0  , \\
  & \bm{y} - \bm{z} = 0  , \\
  & g(\bm{x}, \bm{u}) = 0  , \\
  & h(\bm{x}, \bm{u})  + \bm{s} = 0,  \\
  \label{eq:kkt:cond1}
  & X \bm{v} = 0, \;\bm{x},\bm{v} \geq 0, \\
  \label{eq:kkt:cond2}
  & U \bm{w}= 0, \; \bm{u},\bm{w} \geq 0, \\
  \label{eq:kkt:cond3}
  & S\bm{z}= 0, \; \bm{s},\bm{z} \geq 0.
\end{align}
\end{subequations}

\paragraph{Augmented KKT system.}
The interior-point method replaces the complementarity conditions
by using a homotopy approach. In particular, with a fixed barrier $\mu > 0$
the complementary conditions \eqref{eq:kkt:cond1}-\eqref{eq:kkt:cond3} give
$X \bm{v} = \mu \bm{e}_{n_x}, U \bm{w} = \mu \bm{e}_{n_u}, S \bm{z} = \mu \bm{e}_{m}$
\cite{nocedal_numerical_2006}. Here $\bm{e}_n$ is the vector of all ones of dimension $n$.
By linearizing~\eqref{eq:kkt:fullspace}, we obtain
the following (nonsymmetric) linear system, which is used for the step computation within interior-point iterations.
\begin{subequations}
\begin{equation}
  \label{eq:kktmatrix:original}
  \begin{bmatrix}
    W_{uu} & W_{ux} & 0  & G_u^\top & A_u^\top & -I & 0  & 0  \\
    W_{xu} & W_{xx} & 0  & G_x^\top & A_x^\top & 0  & -I & 0  \\
    0      & 0      & 0  & 0        & -I       & 0  & 0  & -I \\
    G_u    & G_x    & 0  & 0        & 0        & 0  & 0  & 0  \\
    A_u    & A_x    & -I & 0        & 0        & 0  & 0  & 0  \\
    0      & V      & 0  & 0        & 0        & X  & 0  & 0  \\
    W      & 0      & 0  & 0        & 0        & 0  & U  & 0  \\
    0      & 0      & Z  & 0        & 0        & 0  & 0  & S
  \end{bmatrix}
  \begin{bmatrix}
    \bm{p}_u \\
    \bm{p}_x \\
    \bm{p}_s \\
    \bm{p}_\lambda \\
    \bm{p}_y \\
    \bm{p}_{v} \\
    \bm{p}_{w} \\
    \bm{p}_{z}
  \end{bmatrix}
  = -
\begin{bmatrix}
    \bm{g}_u + G_u^\top \bm{\lambda} + A_u^\top \bm{y}  - \bm{w}  \\
    \bm{g}_x + G_x^\top \bm{\lambda} + A_x^\top \bm{y} - \bm{v} \\
    \bm{y} - \bm{z} \\
    g(\bm{x}, \bm{u}) \\
    h(\bm{x}, \bm{u}) + \bm{s}  \\
    X \bm{v} - \mu \bm{e}_{n_x} \\
    U \bm{w} - \mu \bm{e}_{n_u} \\
    S \bm{z} - \mu \bm{e}_m
\end{bmatrix}
\end{equation}
In IPM, it is standard to eliminate the last three block rows
from~\eqref{eq:kktmatrix:original}
(associated with $(\bm{p}_v, \bm{p}_w, \bm{p}_z)$).
The elimination yields the following reduced, symmetric linear system.
\begin{equation}
  \label{eq:kktmatrix:normal}
  \underbrace{
  \begin{bmatrix}
    W_{uu} + \Sigma_u & W_{ux} & 0 & G_u^\top & A_u^\top \\
    W_{xu} & W_{xx} + \Sigma_x & 0 & G_x^\top & A_x^\top \\
    0 & 0 & \Sigma_s & 0 & -I \\
    G_u & G_x & 0 & 0 & 0 \\
    A_u & A_x & -I & 0 & 0
  \end{bmatrix}
}_{K_{aug}}
  \begin{bmatrix}
    \bm{p}_u \\
    \bm{p}_x \\
    \bm{p}_s \\
    \bm{p}_\lambda \\
    \bm{p}_y
  \end{bmatrix}
  = -
  \begin{bmatrix}
    \bm{g}_u + G_u^\top \bm{\lambda} + A_u^\top \bm{y} - \mu U^{-1} \bm{e}_{n_u}  \\
    \bm{g}_x + G_x^\top \bm{\lambda} + A_x^\top \bm{y} - \mu X^{-1} \bm{e}_{n_x} \\
    \bm{y} - \mu S^{-1} \bm{e}_m  \\
    g(\bm{x}, \bm{u}) \\
    h(\bm{x}, \bm{u}) + \bm{s}
  \end{bmatrix}
\end{equation}
\end{subequations}
Here, $\Sigma_u := U^{-1} W$ , $\Sigma_x := X^{-1} V $, and $\Sigma_s:= S^{-1} Z$.
The matrix $K_{aug}$ is
sparse and symmetric indefinite and is typically factorized by a
direct linear solver at each iteration of the interior-point
algorithm. As  discussed before, however, this form is not
amenable to the GPU, motivating us to reduce further
the system~\eqref{eq:kktmatrix:normal}.

\paragraph{Reduced KKT system.}
The reduction strategy adopted in the linearize-then-reduce method
exploits the invertibility of $G_x$. In many applications, such as optimal
power flow, optimal control, and PDE-constrained optimization, the state
variables are completely determined by the control variables. This situation leads to
the structure in the Jacobian block, and in many applications $G_x$ block
is  invertible. In the next theorem, we show how the invertibility of
$G_x$ allows the reduction of the KKT system~\eqref{eq:kktmatrix:normal}.
For ease of  notation,  we denote by $\bm{r} :=
(\bm{r}_1,\bm{r}_2,\bm{r}_3,\bm{r}_4,\bm{r}_5)$ the right-hand side
vector in \eqref{eq:kktmatrix:normal}.
\begin{theorem}{}
  If at $(\bm{x}, \bm{u}) \in \mathbb{R}^{n_u} \times \mathbb{R}^{n_x}$ the Jacobian $G_x$ is invertible, then
  we define the
  reduced Hessian $\hat{W}_{uu}$ and the reduced Jacobian $\hat{A}_u$
  as
  \begin{align*}
    & \hat{W}_{uu} = W_{uu}  - W_{ux} G_x^{-1} G_u - G_u^\top G_x^{-\top} W_{xu} + G_u^\top G_x^{-\top} (W_{xx} + \Sigma_x) G_x^{-1}G_u \; , \\
    & \hat{A}_u = A_u - A_x G_x^{-1}G_u \; .
  \end{align*}
  Then, the augmented KKT system~\eqref{eq:kktmatrix:normal}
  is equivalent to
  \begin{subequations}
  \begin{equation}
    \label{eq:kktmatrix:reduced}
    \underbrace{
    \begin{bmatrix}
      \hat{W}_{uu} + \Sigma_u& 0 & \hat{A}_u^\top \\
      0 & \Sigma_s & -I \\
      \hat{A}_u & - I & 0
    \end{bmatrix}
  }_{K_{red}}
    \begin{bmatrix}
      \bm{p}_u \\
      \bm{p}_s \\
      \bm{p}_y
    \end{bmatrix}
    = -
    \begin{bmatrix}
      \hat{\bm{r}_1} \\ \hat{\bm{r}_2} \\ \hat{\bm{r}_3}
    \end{bmatrix} \; ,
  \end{equation}
  where
  \begin{equation*}
    \left\{
      \begin{aligned}
        & \hat{\bm{r}}_1 = \bm{r}_1 - G_u^\top G_x^{-\top} \bm{r}_2
        - \big(W_{ux} - G_u^\top G_x^{-\top}(W_{xx}+\Sigma_x)\big) G_x^{-1} \bm{r}_4  \; ,\\
        & \hat{\bm{r}}_2 = \bm{r}_3 \; ,\\
        & \hat{\bm{r}}_3 = \bm{r}_5 - A_x G_x^{-1} \bm{r}_4 \; .
      \end{aligned}
    \right.
  \end{equation*}
  Further, the state and adjoint descent directions can be recovered as
  \begin{equation}
    \begin{aligned}
    \bm{p}_x &= - G_x^{-1}\big(\bm{r}_4 + G_u \bm{p}_u\big) \; , \\
    \bm{p}_\lambda &= -G_x^{-\top} \Big(\bm{r}_2 + A_x^\top \bm{p}_y + W_{xu} \bm{p}_u + (W_{xx} + \Sigma_x) \bm{p}_x  \Big) \; .
    \end{aligned}
  \end{equation}
  \end{subequations}
\end{theorem}
{\it Proof.}
Using the fourth block of rows in \eqref{eq:kktmatrix:original},
we can remove the variable $\bm{p}_x$ in the system~\eqref{eq:kktmatrix:reduced}
using $G_x^{-1}$ as a pivot. We get
\begin{equation*}
  \label{eq:kktmatrix:step1}
  \begin{bmatrix}
    W_{uu}+\Sigma_u - W_{ux}G_x^{-1}G_u & 0 & G_u^\top & A_u^\top \\
    W_{xu} - (W_{xx}+\Sigma_x)G_x^{-1}G_u & 0 & G_x^\top & A_x^\top \\
    0 &  \Sigma_s & 0 & -I \\
    A_u - A_x G_x^{-1} G_u & -I & 0 & 0
  \end{bmatrix}
  \begin{bmatrix}
    \bm{p}_u \\
    \bm{p}_s \\
    \bm{p}_\lambda \\
    \bm{p}_y
  \end{bmatrix}
  = -
  \begin{bmatrix}
    \bm{r}_1 - W_{ux}G_x^{-1} \bm{r}_4 \\
    \bm{r}_2 - (W_{xx}+\Sigma_x) G_x^{-1} \bm{r}_4\\
    \bm{r}_3 \\
    \bm{r}_5 - A_x G_x^{-1} \bm{r}_4
  \end{bmatrix}
  \; .
\end{equation*}
The descent direction w.r.t. $\bm{x}$ can be recovered with an additional linear solve:
$\bm{p}_x = - G_x^{-1}\big(\bm{r}_4 + G_u \bm{p}_u\big)$.
In the system~\eqref{eq:kktmatrix:step1} we can eliminate the
third block of columns (associated with $\bm{p}_\lambda$), with $G_x^{-\top}$ as a pivot.
We recover the linear system in~\eqref{eq:kktmatrix:reduced}.
The descent direction $\bm{p}_\lambda$ now satisfies
\begin{equation*}
  \begin{aligned}
    \bm{p}_\lambda &= -G_x^{-\top}\Big(
    \bm{r}_2 - (W_{xx}+\Sigma_x)G_x^{-1} \bm{r}_4
  + A_x^\top \bm{p}_y + (W_{xu} - (W_{xx}+\Sigma_x) G_x^{-1} G_u)\bm{p}_u \Big) \\
                   &= -G_x^{-\top} \Big(\bm{r}_2 + A_x^\top \bm{p}_y + W_{xu} \bm{p}_u + (W_{xx} + \Sigma_x) \bm{p}_x  \Big)
  \end{aligned}
\end{equation*}
by rearranging the terms,
thus completing the proof.
\qed

\paragraph{Condensed KKT system.}
The left-hand side matrix $K_{red}$ in~\eqref{eq:kktmatrix:reduced}
has a size $(2n_u + m)\times (2n_u + m)$, which can be prohibitively large
if the number of constraints $m$ is substantial. Fortunately,
\eqref{eq:kktmatrix:reduced} can be condensed further,
down to a system with size $n_u \times n_u$.
This additional reduction requires the invertibility of $\Sigma_s$, which is always
the case for interior point algorithms.
\begin{theorem}{}
  \label{thm:kkt:condensed}
  We suppose that $\Sigma_s$ is nonsingular.
  Then the linear system~\eqref{eq:kktmatrix:reduced} is equivalent to
  \begin{subequations}
  \begin{equation}
    \label{eq:kktmatrix:condensed}
    \underbrace{
    \Big(\hat{W}_{uu} + \hat{A}_u^\top \Sigma_s \hat{A}_u \Big)
  }_{K_{cond}}
    \bm{p}_u
    = - \big(\hat{\bm{r}}_1 + \hat{A}_u^\top \Sigma_s \hat{\bm{r}}_3 - \hat{A}_u^\top \hat{\bm{r}}_2 \big)
    \; ,
  \end{equation}
  the slack and multiplier descent directions being recovered as
  \begin{equation}
    \bm{p}_y = \Sigma_s \big(\hat{A}_u \bm{p}_u + \hat{\bm{r}}_3 - \Sigma_s^{-1} \hat{\bm{r}}_2 \big) \; , \quad
    \bm{p}_s = \Sigma_s^{-1} (\bm{p}_y - \hat{\bm{r}}_2 )\;.
  \end{equation}
  \end{subequations}
\end{theorem}
{\it Proof.}
  In~\eqref{eq:kktmatrix:reduced} we eliminate the second block of columns (associated
  with the slack descent $\bm{p}_s$) to get
  \begin{equation}
    \label{eq:kktmatrix:tmp2}
    \begin{bmatrix}
      \hat{W}_{uu} &  \hat{A}_u^\top \\
      \hat{A}_u & - \Sigma_s^{-1}
    \end{bmatrix}
    \begin{bmatrix}
      \bm{p}_u \\
      \bm{p}_y
    \end{bmatrix}
    = -
    \begin{bmatrix}
      \hat{\bm{r}}_1 \\
      \hat{\bm{r}}_3 - \Sigma_s^{-1} \hat{\bm{r}}_2
    \end{bmatrix} \; .
  \end{equation}
  The slack descent direction is recovered as
  $
    \bm{p}_s = \Sigma_s^{-1} (\bm{p}_y - \hat{\bm{r}}_2 )\; .
  $
  In the last block of rows we can reuse $\Sigma_s^{-1}$ as a pivot to simplify
  further~\eqref{eq:kktmatrix:tmp2}:
  \begin{equation}
    \Big(\hat{W}_{uu} + \hat{A}_u^\top \Sigma_s \hat{A}_u \Big) \bm{p}_u =
    - (\hat{\bm{r}}_1 + \hat{A}_u^\top \Sigma_s \hat{\bm{r}}_3 - \hat{A}_u^\top \hat{\bm{r}}_2 )
  \end{equation}
  with
  $\bm{p}_y = \Sigma_s \big(\hat{A}_u \bm{p}_u + \hat{\bm{r}}_3 - \Sigma_s^{-1} \hat{\bm{r}}_2 \big)$.
\qed

\paragraph{Discussion.}
The final matrix $K_{cond}$ has a size $n_u \times n_u$ and
is dense, meaning that it can be factorized efficiently by any LAPACK library.
We note that the reduction has proceeded in two steps: first, we have reduced
the system by eliminating the state variables, and then we have condensed it by eliminating the slacks, giving the order $K_{aug} \to K_{red} \to
K_{cond}$. We have opted for this order to simplify the comparison
with the reduce-then-linearize approach presented in the next
section~\ref{sec:reduced:fripm}. In practice, however, it is more convenient
to first condense the linear system and then reduce it: $K_{aug} \to K_{cond} \to K_{red}$.
This equivalent approach avoids the allocation of the dense reduced Jacobian $\hat{A}_u$, which has a size $m \times n_u$ and is expensive to store in memory.
The reduction $K_{aug} \to K_{cond} \to K_{red}$ is illustrated
in Figure~\ref{fig:reduced:kktreduction}.

\begin{figure}[!ht]
  \centering
  \includegraphics[width=12cm]{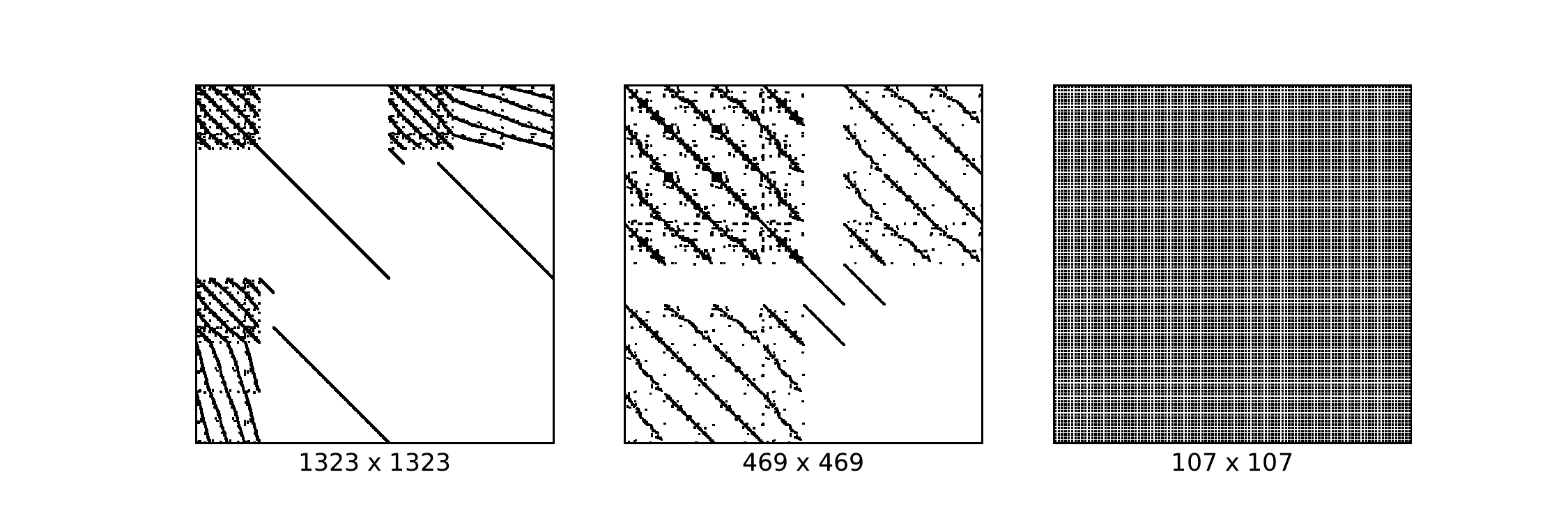}
  \caption{Successive reductions $K_{aug} \to K_{cond} \to K_{red}$ associated with {\tt 118ieee}.}
  \label{fig:reduced:kktreduction}
\end{figure}

We establish the last condition to guarantee that we compute
a descent direction at each iteration.
For any symmetric matrix $M \in \mathbb{R}^{n \times n}$,
we note its inertia $\text{I}(M) = (n_+, n_-, n_0)$ as
the numbers of positive, negative, and zero eigenvalues.
\begin{theorem}
  \label{thm:reduction:descentdirection}
  The step $\bm{p}$ in \eqref{eq:kktmatrix:normal} is a descent direction if
  \begin{itemize}
    \item $\text{I}(K_{aug}) = (n_x+n_u+m, n_x + m, 0)$, equivalent to
    \item $\text{I}(K_{red}) = (n_u+m, m, 0)$, equivalent to
    \item $\text{I}(K_{cond}) = (n_u, 0, 0)$ (the condensed matrix $K_{cond}$ is
      positive definite).
  \end{itemize}
\end{theorem}
The equivalence of the three conditions can be verified via the Haynsworth inertia additivity formula.

\subsubsection{Linearize-then-reduce algorithm (LinRed IPM)}
Now that the different reductions have been introduced,
we are able to present the linearize-then-reduce ({\tt LinRed}) algorithm
in Algorithm~\ref{algo:ripm}, following \cite{cervantes2000reduced}.
The {\bf reduction} step is in itself an expensive operation,
as we will explain in \S\ref{sec:streamline:reduction}. The other bottlenecks are the factorization of the
dense condensed matrix $K_{cond}$ (which amounts to a Cholesky
factorization if the matrix is positive definite)
and the factorization of the sparse Jacobian $G_x$.

\begin{algorithm}
  \KwData{Initial primal variables $(\bm x_0, \bm u_0, \bm s_0)$ and dual variables $(\bm \lambda_0, \bm y_0)$}
\For{$k = 0,...$}{
  Evaluate the derivatives $W, G, A$ at $(\bm x_k, \bm u_k)$ \;
  {\bf Reduction:} Condense the KKT system \eqref{eq:kktmatrix:normal} in $K_{cond}$ \;
  {\bf Control step:} Factorize~\eqref{eq:kktmatrix:condensed} and solve the system to find $\bm p_u$ \;
  {\bf Dual step:} $\bm{p}_y = \Sigma_s (\hat{A}_u \bm{p}_u + \hat{\bm{r}}_3 - \Sigma_s^{-1} \hat{\bm{r}}_2 )$\;
  {\bf Slack step:} $\bm{p}_s = \Sigma_s^{-1} (\bm{p}_y - \hat{\bm{r}}_2 )$\;
  {\bf State step:} $\bm{p}_x = - G_x^{-1}\big(\bm{r}_4 + G_u \bm{p}_u\big)$ \;
  {\bf Adjoint step:} $\bm{p}_\lambda =  -G_x^{-\top} (\bm{r}_2 + A_x^\top \bm{p}_y + W_{xu} \bm{p}_u + (W_{xx} + \Sigma_x) \bm{p}_x  )$ \;
  {\bf Line search:} Update primal-dual direction $(\bm u_{k+1}, \bm{x}_{k+1}, \bm{s}_{k+1},
  \bm{\lambda}_{k+1}, \bm{y}_{k+1})$
  using a filter line-search along direction
  $(\bm{p}_u,\bm{p}_x,\bm{p}_s,\bm{p}_\lambda,\bm{p}_y)$\;
}
 \caption{Linearize-then-reduce algorithm}
 \label{algo:ripm}
\end{algorithm}

\subsection{Reduce-then-linearize}
\label{sec:reduced:fripm}
We now focus on our second reduction scheme, operating directly
at the level of the nonlinear problem~\eqref{eq:originalproblem}.
This method can be interpreted as an interior-point
alternative of the reduced-gradient algorithm of Dommel and Tinney~\cite{dommel_optimal_1968}
and was explored recently in \cite{kardos2020reduced}.

\subsubsection{Nonlinear reduction}
\paragraph{Nonlinear projection.}
Instead of operating the reduction in the linearized KKT system as
in \S\ref{sec:reduced:ripm},
Dommel and Tinney's method uses the implicit function theorem
to remove the state variables from the problem.
\begin{theorem}{(Implicit function theorem).}
  Let $g: \mathbb{R}^{n_x} \times \mathbb{R}^{n_u} \to \mathbb{R}^{n_x}$
  a continuously differentiable function, and let $(\bm{x}, \bm{u})
  \in \mathbb{R}^{n_x} \times \mathbb{R}^{n_u}$
  such that $g(\bm x, \bm u) = 0$. If the Jacobian $G_x$
  is \emph{invertible}, then there exist an open set
  $U \subset \mathbb{R}^{n_u}$ (with $\bm{u} \in U$)
  and a unique differentiable function $\underline{x}: U \to \mathbb{R}^{n_x}$
  such that $g(\underline{x}(\bm u), \bm u) = 0$ for all $\bm u \in U$.
\end{theorem}
The implicit function theorem gives, under certain assumptions,
the existence of a local differentiable function $\underline{x}: U \to \mathbb{R}^{n_x}$
attached to a given control~$\bm u$.
If we assume that on the feasible domain $g(\bm x, \bm u) = 0$
is solvable and the Jacobian is invertible everywhere,
we can derive the reduced-space problem as
\begin{equation}
  \label{eq:reduced:problem}
  \min_{\bm{u}} \quad  f(\underline{x}(\bm u), \bm u) \quad
  \text{subject to} \quad
  \left\{
  \begin{aligned}
     \quad  & \bm{u} \geq 0 \, , \quad \underline{x}(\bm u) \geq 0 \\
                             &  h(\underline{x}(\bm u), \bm u) \leq 0 \, .
  \end{aligned}
  \right.
\end{equation}
In contrast to~\eqref{eq:originalproblem},
the problem~\eqref{eq:reduced:problem} optimizes only
with relation to the control~$\bm u$, the state being defined
\emph{implicitly} via the local functionals $\underline{x}$.
By definition $g(\underline{x}(\bm u), \bm u) = 0$, the
state equation is automatically satisfied in the reduced-space.
However, the reduced-space problem is tied
to the
assumptions of the implicit function theorem:
in some applications,
finding a control $\bm u$ invertible w.r.t. the
state equations~\eqref{eq:stateequation} can be challenging.

\paragraph{Reduced derivatives.}
We define the reduced objective and the reduced constraints as
\begin{equation}
  f_r(\bm u) := f(\underline{x}(\bm u), \bm u) \, , \quad
  h_r(\bm u) :=
  \begin{bmatrix}
  h(\underline{x}(\bm u), \bm u)  \\
  - \underline{x}(\bm u)
  \end{bmatrix}
  \, ,
\end{equation}
and we note the reduced Lagrangian
$L_r(\bm u, \bm s ; \bm y) := f_r(\bm u) + \bm{y}^\top (h_r(\bm u) + \bm s)
- \bm{w}^\top \bm{u} - \bm{z}^\top \bm{s}$.
By exploiting the implicit function theorem, we can deduce the derivatives in the reduced-space.

\begin{theorem}{(Reduced derivatives~\cite[Chapter 15.]{griewank2008evaluating})}
  Let $\bm u \in \mathbb{R}^{n_u}$ such that the conditions of the
  implicit function theorem hold. Then,
  \begin{subequations}
        \label{eq:reduced:derivatives}
  \begin{itemize}
    \item The functions $f_r$ and $h_r$ are
      continuously differentiable, with
      \begin{equation}
        \nabla f_r(\bm u) = \bm{g}_u - G_u^\top G_x^{-\top} \bm{g}_x  \, , \quad
        \hat{A}_u =  \partial h_r(\bm u) =
        \begin{bmatrix}
        A_u - A_x G_x^{-1} G_u  \\
         G_x^{-1} G_u
        \end{bmatrix}
        \; .
      \end{equation}
    \item The Hessian of the reduced Lagrangian $L_r$ satisfies
      \begin{equation}
        \hat{W}_{uu} = W_{uu} - W_{ux} G_x^{-1} G_u - G_u^\top G_x^{-\top} W_{xu} + G_u^\top G_x^{-\top} W_{xx} G_x^{-1}G_u \; .
      \end{equation}
  \end{itemize}
  \end{subequations}
\end{theorem}

\paragraph{Reduced-space KKT system.}
The KKT conditions associated with the reduced problem \eqref{eq:reduced:problem} are
\begin{subequations}
  \label{eq:kkt:reduced}
\begin{align}
  & \nabla_{u} f_r + \hat{A}_u^\top \bm{y} - \bm{w} = 0 , \\
  & \bm{y} - \bm{z} = 0 , \\
  & h_r(\bm{u})  + \bm{s} = 0 , \\
  & U\bm{w}= 0 , \;\bm{u},\bm{w} \geq 0, \\
  & S\bm{z}= 0 , \;\bm{s},\bm{z} \geq 0 ,
\end{align}
\end{subequations}
translating, at each iteration of the IPM algorithm, to
the following augmented KKT system:
\begin{equation}
  \label{eq:reduced:kktsystem}
  \begin{bmatrix}
    \hat{W}_{uu} + \Sigma_u & \phantom{00} 0\phantom{00}  & \hat{A}_u^\top \\
    0 & \Sigma_s & -I \\
    \hat{A}_u & -I & 0
  \end{bmatrix}
  \begin{bmatrix}
    \bm{p}_u \\ \bm{p}_s \\ \bm{p}_y
  \end{bmatrix}
  = -
  \begin{bmatrix}
    \nabla_u L_r - \mu U^{-1} \bm{e}_{n_u} \\
    \bm{y} - \mu S^{-1} \bm{e}_{m} \\
    h_r(\bm{u}) + \bm{s}
  \end{bmatrix}
  \, .
\end{equation}
We note that we can apply Theorem~\ref{thm:kkt:condensed}
to get a condensed form of the KKT system \eqref{eq:reduced:kktsystem}.
The KKT system~\eqref{eq:reduced:kktsystem} has a  structure
similar to that of \eqref{eq:kktmatrix:reduced} but  with minor differences.
(i) The state $\bm{x}$ and the adjoint $\bm{\lambda}$ are
updated independently of \eqref{eq:reduced:kktsystem}, respectively
by solving the state equation~\eqref{eq:stateequation} and by solving a linear system.
(ii) The bounds $\bm{x} \geq 0$ are incorporated inside the
function $h_r$, whereas \eqref{eq:kktmatrix:reduced} handles them
explicitly in the KKT system (leading to an additional term $\Sigma_x$
in the reduced Hessian $\hat{W}_{uu}$).
(iii) The right-hand side in \eqref{eq:kktmatrix:reduced} incorporates
additional second-order terms that do not appear in
\eqref{eq:reduced:kktsystem}.
Indeed, as here with $\bm{r}_4 = g(\bm x, \bm u) = 0$, all the
terms associated with $\bm{r}_4$ disappeared in the right-hand side
of \eqref{eq:reduced:kktsystem}, including the second-order terms.

\subsubsection{Reduce-then-linearize algorithm (RedLin IPM)}
At each iteration $k$,
the reduce-then-linearize ({\tt RedLin}) algorithm proceeds in two steps.
First, given a new control $\bm u_k$, the algorithm finds a state
$\bm x_k$ satisfying $g(\bm x_k, \bm u_k) = 0$ by using
a nonlinear solver. Then, once $\bm x_k$ has been computed, the reduced
derivatives are updated, and the condensed form of the system~\eqref{eq:reduced:kktsystem}
is solved to compute the next iterate. This process is summarized in Algorithm~\ref{algo:fripm}.
We note that compared with Algorithm~\ref{algo:ripm}, the
{\bf state step} and the {\bf adjoint step} are computed before
solving the KKT system, since we first have to solve~\eqref{eq:stateequation}.
In addition, the algorithm gives no guarantee that at iteration
$k$ there exists a state $\bm{x}_k$ such that $g(\bm x_k, \bm u_k) = 0$,
which can be problematic on certain applications. On the other hand, even if
interrupted early, the algorithm produces state variables that are feasible
for the state equation~\eqref{eq:stateequation}, which may be important in real-time applications.

\begin{algorithm}
  \KwData{Initial primal variable $(\bm u_0, \bm s_0)$ and dual variable $(\bm y_0)$}
\For{$k = 0,...$}{
  {\bf Projection:} Find $\bm x_k$ satisfying $g(\bm x_k, \bm u_k) = 0$ \;
  {\bf Adjoint step:} Solve $\bm{\lambda} = -G_x^{-\top} \nabla_x L$ \;
  {\bf Reduction:} Condense the KKT system \eqref{eq:reduced:kktsystem} in $K_{cond}$ \;
  {\bf Control step:} Factorize~$K_{cond}$ and solve the system to find $\bm p_u$ \;
  {\bf Dual step:} $\bm{p}_y = \Sigma_s (\hat{A}_u \bm{p}_u + \hat{\bm{r}}_3 - \Sigma_s^{-1} \hat{\bm{r}}_2 )$\;
  {\bf Slack step:} $\bm{p}_s = \Sigma_s^{-1} (\bm{p}_y - \hat{\bm{r}}_2 )$\;
  {\bf Line search:} Update primal-dual direction $(\bm u_{k+1}, \bm{s}_{k+1},
  \bm{y}_{k+1})$
  using a filter line-search along direction
  $(\bm{p}_u,\bm{p}_s,\bm{p}_y)$\;
}
 \caption{Reduce-then-linearize algorithm}
 \label{algo:fripm}
\end{algorithm}

\section{Implementation of reduced IPM on GPU}
\label{sec:streamline}

In this section we present a GPU implementation of the reduced-space algorithm.
To avoid expensive data transfers between the host and the device,
we have designed our implementation to run as much as possible on the GPU, comprising
(i) the evaluation of the callbacks, (ii) the reduction algorithm,
and (iii) the dense factorization of the condensed KKT system.
At the end, only the interior-point routines (line search, second-order correction, etc.)
run on the host.
In particular, our method does not require transferring dense matrices between device
and host, and most of the operations performed on the host are simple scalar operations.

\subsection{GPU operations}
The key idea is to exploit efficient computation kernels on the GPU
to implement the reduction algorithm.
To the extent possible, we avoid writing custom kernels and rely instead on the BLAS and
LAPACK operations, as provided by the vendor library (CUDA in our case).

We list below the specific kernels we are targeting ({\tt Sp} stands for \emph{sparse}, {\tt Dn} for \emph{dense}).
\begin{itemize}
  \item {\tt SpMV/SpMM} \emph{(sparse matrix-vector product/sparse matrix-dense matrix product)}. On the GPU, sparse matrices are stored in condensed sparse row
    (CSR) format, allowing  the
    sparse multiplication kernels to be run fully in parallel.
  \item {\tt SpSV/SpSM} \emph{(sparse triangular solve)}. Once an LU factorization is
    computed (for instance with {\tt SpRF}), a sparse matrix $A$ is decomposed
    as $P A Q = LU$, with $P$ and $Q$ two permutation matrices, $L$ and $U$
    being respectively a lower and an upper triangular matrix. The routine {\tt SpSV}
    solves the triangular systems $L^{-1}\bm{b}$ and $U^{-1} \bm{b}$ efficiently.
    The extension to multiple right-hand sides is directly provided by
    the {\tt SpSM} kernel.
  \item {\tt SpRF} \emph{(sparse LU refactorization)}. GPUs are notoriously
    inefficient at factorizing a sparse matrix $M$. If the sparse matrix
    $M$ has always the same sparsity pattern, however, we can compute the initial
    factorization on the CPU and move the factors back to the GPU.
    Then, if the nonzero coefficients of the matrix changes,
    the matrix can be refactorized entirely on the GPU, by updating
    directly the $L$ and $U$ factors.
\end{itemize}

\subsection{Porting the callbacks to the GPU}
\label{sec:streamline:autodiff}
In nonlinear programming, the evaluation of the callbacks is often
one of the most time-consuming parts,
and the OPF problem is not immune to this issue.
Most of the time, the derivatives of the OPF model are provided
explicitly~\cite{zimmerman2010matpower}, evaluated by using automatic differentiation~\cite{dunning2017jump,fourer1990modeling}
or  symbolic differentiation~\cite{hijazi2018gravity}.
Here we have chosen to stick with automatic differentiation.
We have streamlined on the GPU the evaluation of the objective
and of the constraints
by adopting the vectorized model proposed in~\cite{lee2020feasible}.
This model factorizes all the nonlinearities inside a single basis function,
which can be evaluated in parallel inside a single GPU kernel.
In addition, we have designed our implementation to be differentiable with {\tt ForwardDiff.jl}~\cite{revels2016forward}.
In total, each iteration of the reduced IPM algorithm requires one
evaluation of the objective's gradient (=one reverse pass),
one evaluation of the
Jacobian of the constraints (=one forward pass), and
one evaluation of the Hessian of the Lagrangian (=one forward-over-reverse pass).

The OPF comes with two decisive advantages. (i) Its structure is supersparse,
rendering all the {\tt SpMV} operations efficient (we have at most
a dozen of nonzeroes on each rows of the sparse matrices).
(ii) The sparsity of the Jacobian $G_x$ is fixed (it is associated with
the structure of the power network), allowing for efficient refactorization
with {\tt SpRF}.

\subsection{Porting the reduction algorithm to the GPU}
\label{sec:streamline:reduction}
Once the derivatives are evaluated in the full space,
it remains to build the condensed KKT system~\eqref{eq:kktmatrix:condensed}.
The algorithm has to be repeated at each iteration of the IPM algorithm,
and its performance is critical.
If we denote $K = W + A^\top \Sigma_s A$, then we note that
\begin{equation}
  \label{eq:algo:reduction}
   \hat{W}_{uu} + \hat{A}_u^\top \Sigma_s \hat{A}_u =
 \begin{bmatrix}
   I \\ - G_x^{-1} G_u
 \end{bmatrix}^\top
 \begin{bmatrix}
   K_{uu} & K_{ux} \\
   K_{xu} & K_{xx}
 \end{bmatrix}
 \begin{bmatrix}
   I \\ - G_x^{-1} G_u
 \end{bmatrix}
 \equiv \hat{K}_{uu}
 \; .
\end{equation}
Evaluating \eqref{eq:algo:reduction} requires three different
operations: (i) factorizing the Jacobian $G_x$ ({\tt SpRF}), (ii)
triangular solves $G_x^{-1} b$ ({\tt SpSV/SpSM}), and (iii)
sparse-matrix matrix multiplications with $K$ ({\tt SpMM}).
The order in which the operations are performed is
important, affecting the complexity of the reduction algorithm.

The naive idea is to evaluate first the sensitivity matrix $S=-G_x^{-1}G_u$,
in order to reduce the total number of linear solves to $n_u$. However,
the matrix $S$ is dense, with size $n_x \times n_u$, and the
cross-product $S^\top K_{xx} S$ requires storing another intermediate dense matrix with size $n_x \times n_u$ to evaluate the dense product $S^{\top} (K_{xx} S)$ (which is itself slow when $n_x$ is large). This renders the algorithm
not tractable on the largest instances.

Hence, we avoid computing the full sensitivity matrix $S$ and
rely instead on a batched variant of the adjoint-adjoint algorithm~\cite{pacaud2022batched}.
First,
we  compute an LU factorization of $G_x$, as $P G_x Q = L \, U$, with $P$ and $Q$
two permutation matrices and $L$ and $U$ being respectively a lower
and an upper triangular matrix
(using {\tt SpRF}, the factorization can be updated entirely on the GPU if the sparsity
pattern of $G_x$ is the same along the iterations). Once the factorization is computed,
solving the linear solve $G_x^{-1} \bm{b}$ translates to
2 {\tt SpMV} and 2 {\tt SpSV} routines, as $G_x^{-1} \bm{b} = QU^{-1} L^{-1} P \bm{b}$.

Second, we build the sparse matrix $K
\in \mathbb{R}^{(n_x+n_u) \times (n_x+n_u)}$ (nontrivial but doable in one sparse addition and one
sparse-sparse multiplication {\tt SpGEMM}).
Then, for a batch size $N$, the algorithm takes as input a dense matrix $V \in \mathbb{R}^{n_u \times N}$
and evaluates the Hessian-matrix product $\hat{K}_{uu} V$
with three successive operations.
\begin{enumerate}
  \item Solve $Z = -G_x^{-1} (G_u \, V)$.
    \hfill (3 {\tt SpMM}, 2 {\tt SpSM})
  \item Evaluate
    $\begin{bmatrix}
        H_u \\ H_x
      \end{bmatrix}
      =
      \begin{bmatrix}
        K_{uu} & K_{ux} \\
        K_{xu} & K_{xx}
      \end{bmatrix}
      \begin{bmatrix}
        V \\ Z
      \end{bmatrix}
    $.
    \hfill (1 {\tt SpMM})
  \item Solve $\Psi = G_x^{-\top} H_x$ and
    get $\hat{K}_{uu} = H_u - G_u \Psi$.
    \hfill (3 {\tt SpMM}, 2 {\tt SpSM})
\end{enumerate}
One Hessian-matrix product $\hat{K}_{uu}V$
requires $2 N$ linear solves (streamlined
in four {\tt SpSM} operations),
giving a total of 7 {\tt SpMM} and
4 {\tt SpSM} operations.
Evaluating $\hat{K}_{uu}V$ requires
 $3 \times n_x \times N$ storage for
the two intermediates $Z, \Psi \in \mathbb{R}^{n_x \times N}$, as well as
an additional buffer to store the permuted matrix in the LU triangular solves.
Overall, the evaluation of the full reduced matrix $\hat{K}_{uu}$ requires
$\text{div}(n_u, N) + 1$ Hessian-matrix products, giving a complexity
proportional to the number of controls $n_u$.

\section{Numerical results}
\label{sec:results}
In this section we assess the performance of the reduced IPM
algorithm on the various OPF instances presented
in \S\ref{sec:results:instances}.
First, we evaluate in \S\ref{sec:results:reduction} the
scalability of the GPU-accelerated reduction algorithm initially
presented in \S\ref{sec:streamline:reduction}.
Next, we present in \S\ref{sec:results:linred} a detailed
assessment of the reduced IPM algorithm on the GPU,
by comparing its performance with that of a state-of-the-art full-space IPM working on the CPU.
Then, we present in \S\ref{sec:results:redlin} a comparison of
our two reduced-space algorithms: \textit{linearize-then-reduce} ({\tt LinRed})
and \textit{reduce-then-linearize} ({\tt RedLin}).

\subsection{Benchmark instances}
\label{sec:results:instances}
We benchmark the reduced IPM algorithm on the OPF problem.
We select in Table~\ref{tab:test_instances} a subset of the MATPOWER cases
provided in \cite{zimmerman2010matpower}, whose size varies from medium to large scale.
In addition, we add three cases from the PGLIB benchmark~\cite{babaeinejadsarookolaee2019power} with fewer
degrees of freedom (as indicated by the ratio $\frac{n_u}{n_x + n_u}$).
Our algorithms have been implemented
in Julia for portability. All the benchmarks presented have been
generated on our workstation, equipped with an NVIDIA V100 GPU
and using {\tt CUDA 11.4}.
The code, open source, is available on \url{https://github.com/exanauts/Argos.jl/}.

\begin{table}[!ht]
  \centering
  \resizebox{.6\textwidth}{!}{
\begin{tabular}{l|rrr|rrc}
\toprule
Case &     $n_b$ &     $n_\ell$ & $n_g$ &      $n_x$ &     $n_u$  & $n_u / (n_x+n_u)$\\
\midrule
                  118ieee & 118   & 186   & 54   & 181    & 107   & 0.37 \\
                  300ieee & 300   & 411   & 69   & 530    & 137   & 0.21 \\
          ACTIVSg500      & 500   & 597   & 56   & 943    & 111   & 0.11 \\
           1354pegase     & 1,354  & 1,991  & 260  & 2,447   & 519   & 0.17 \\
         ACTIVSg2000      & 2,000  & 3,206  & 432  & 3,607   & 783   & 0.18 \\
           2869pegase     & 2,869  & 4,582  & 510  & 5,227   & 1,019  & 0.16 \\
           9241pegase     & 9,241  & 16,049 & 1,445 & 17,036  & 2,889  & 0.14 \\
          ACTIVSg10k      & 10,000 & 12,706 & 1,937 & 18,544  & 2,909  & 0.14 \\
          13659pegase     & 13,659 & 20,467 & 4,092 & 23,225  & 8,183  & 0.26 \\
          ACTIVSg25k      & 25,000 & 32,230 & 3,779 & 47,246  & 5,505  & 0.10 \\
          ACTIVSg70k      & 70,000 & 88,207 & 8,107 & 134,104 & 11,789 & 0.08 \\
\midrule
   9591\_goc              & 9,591  & 15,915 & 365  & 19,013  & 335   & 0.02 \\
10480\_goc                & 10,480 & 18,559 & 777  & 20,620  & 677   & 0.03 \\
  19402\_goc              & 19,402 & 34,704 & 971  & 38,418  & 769   & 0.02 \\
\bottomrule
\end{tabular}
}
\caption{Case instances obtained from MATPOWER}
\label{tab:test_instances}
\end{table}

\subsection{How far can we parallelize the reduction algorithm on the GPU?}
\label{sec:results:reduction}
We first benchmark the reduction algorithm on the instances
in Table~\ref{tab:test_instances}. For {\tt SpRF}, the reduction algorithm
uses the library {\tt cusolverRF} (with an initial factorization computed by KLU),
the kernels {\tt SpSM} and {\tt SpMM} being provided by {\tt cuSPARSE}.
We show
in Figure~\ref{fig:results:reduction} (a) that {\tt cusolverRF} is able to refactorize efficiently the Jacobian $G_x$
on the GPU (its sparsity pattern is constant, and given by the structure of the underlying network).
In Figure~\ref{fig:results:reduction} (b), we depict the performance of the
reduction algorithm against the batch size $N$. We observe that the greater the size $N$, the better is the performance, until we reach the scalability limit of the GPU.
For instance, on {\tt ACTIVSg70k} we reach the limit when $N=512$, meaning
we cannot parallelize the algorithm further on the GPU. This limits the performance of the reduction
since $\text{div}(11\,789, 512) + 1 = 24$ batched Hessian-matrix products remain
to be computed to evaluate the reduced Hessian of {\tt ACTIVSg70k}.
We note that, overall, it makes no sense to use a batch size greater
than $N=256$.

\begin{figure}[!ht]
  \centering
  \begin{tabular}{cc}
    \includegraphics[width=.4\textwidth]{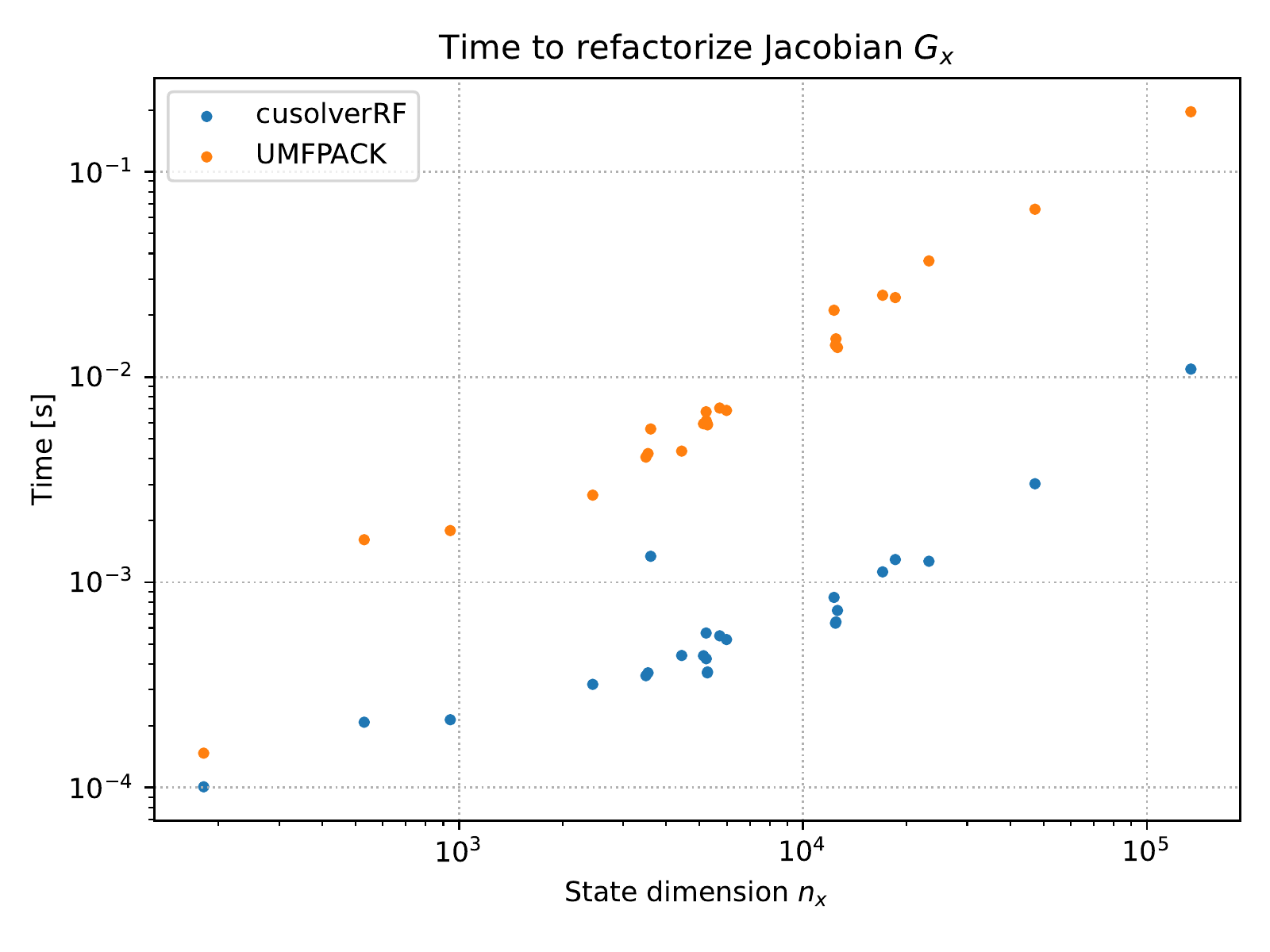} &
    \includegraphics[width=.4\textwidth]{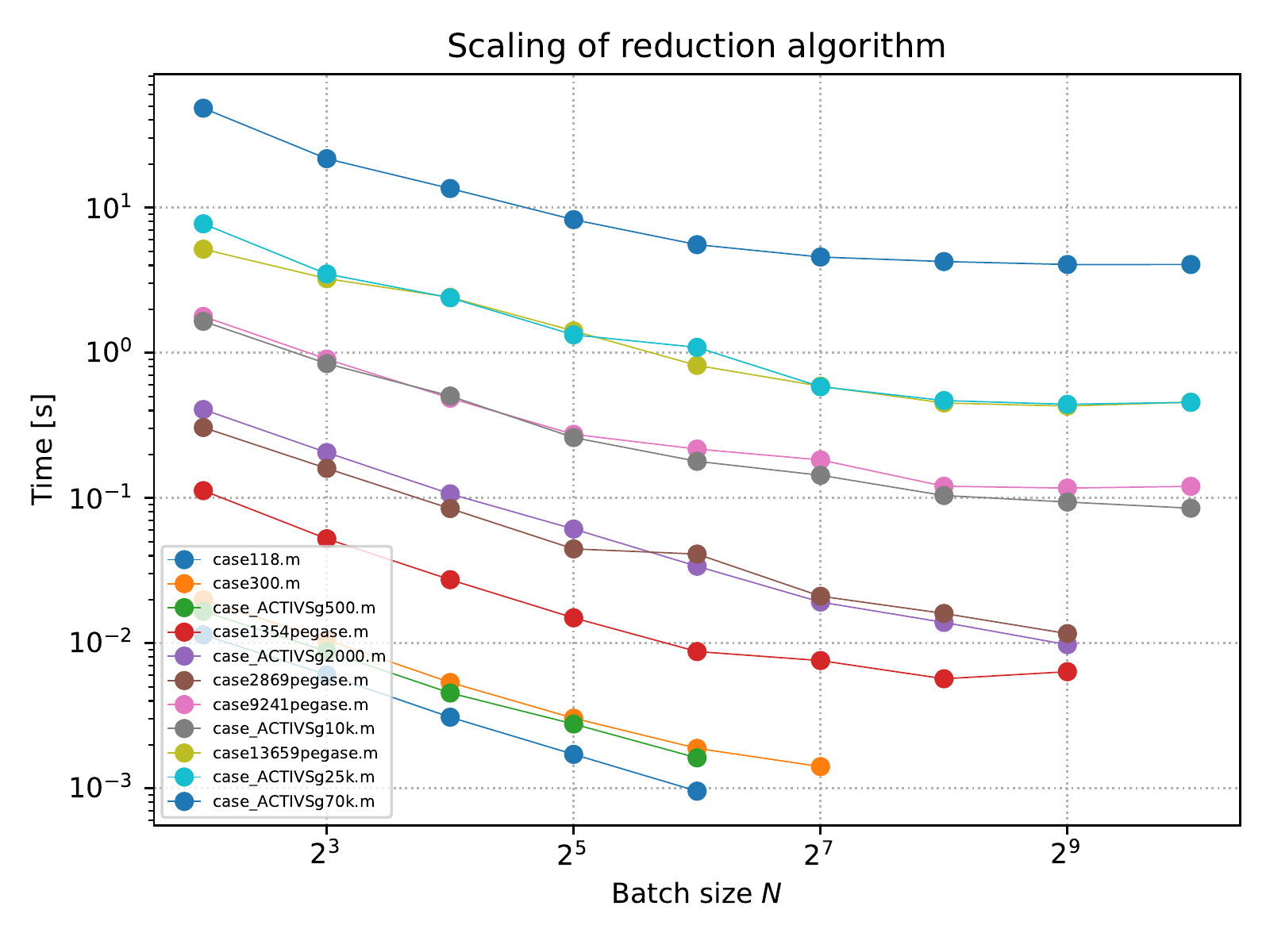} \\
    (a) & (b)
  \end{tabular}
  \caption{Performance of the reduction algorithm: (a) time spent refactorizing
  the Jacobian~$G_x$ in {\tt cusolverRF} and in {\tt UMFPACK}; (b) performance of the reduction algorithm against the batch size $N$.}
  \label{fig:results:reduction}
\end{figure}


\subsection{Linearize-then-reduce or full-space interior-point?}
\label{sec:results:linred}
We have implemented the \emph{linearize-then-reduce} ({\tt LinRed}) inside the MadNLP solver~\cite{shin2021graph}.
As a reference, we benchmark our code against Knitro 13.0~\cite{waltz2006interior}
and MATPOWER~\cite{zimmerman2010matpower}. We note that the performance
is consistent with that reported in the recent benchmark
\cite{kardos2018complete}. Our implementation uses the vectorized OPF model
introduced in \S\ref{sec:streamline:autodiff}, which runs entirely on the GPU
(including the evaluation of the derivatives).
MadNLP runs in inertia-based mode: Theorem~\ref{thm:reduction:descentdirection}
states that the inertia is correct if and only if the reduced matrix $K_{cond}$ is
positive-definite. This fact is exploited in {\tt LinRed}. At
each iteration the algorithm factorizes the matrix $K_{cond}$ with the dense Cholesky
solver shipped with {\tt cusolver}; if the factorization fails, we apply a primal-dual regularization
to $K_{cond}$ until it becomes positive-definite.
In addition, we use the reduction algorithm presented
in \S\ref{sec:streamline:reduction} with a batch size $N=256$.
We compare {\tt LinRed} with a hybrid \emph{full-space IPM} algorithm
using also the GPU-accelerated OPF model but solving the original augmented
system~\eqref{eq:kktmatrix:normal} on the CPU with the linear solver MA27
(hence combining the best of both worlds).
The hybrid \emph{full-space} IPM and {\tt LinRed} both use the same model
and the same derivatives and are equivalent in exact arithmetic unless we run
into a feasibility restoration phase. In that edge case, the algebra of {\tt LinRed}
can be adapted but no longer follows the workflow we presented above.
In that circumstance, {\tt LinRed} applies the dual regularization only on the inequality constraints,
to fit the linear algebra framework we introduced in \S\ref{sec:reduced:ripm}.

Concerning scaling, MadNLP uses the same approach as Ipopt, based on the norm of
the first-order derivatives~\cite{wachter2006implementation}.
The initial primal variables $(\bm{x}_0, \bm{u}_0)$ are specified
inside the MATPOWER file, and the initial dual variables are set to zero: $\bm{y}_0 = 0$.
The algorithm stops when the primal and dual infeasibilities are below $10^{-8}$.
The results of the benchmark are presented in Table~\ref{tab:results:benchmark}.
Here, the dagger sign $\dagger$ indicates that the IPM algorithm runs into a feasibility restoration: this is the case both for {\tt ACTIVSg10k} and for
{\tt ACTIVSg70k} (even Knitro struggles on these two instances, with
numerous conjugate gradient iterations performed).

\begin{table}[!ht]
  \centering
  \resizebox{\textwidth}{!}{
  \begin{tabular}{l|rr|rrr|rrrr}
    \toprule
    & \multicolumn{2}{c}{Knitro} &\multicolumn{3}{c}{Hybrid full-space IPM} &  \multicolumn{4}{c}{\tt LinRed} \\
    \midrule
                  Case  & \#it & Time (s)     & \#it & Time (s)      & MA27 (s) & \#it & Time (s)     & Chol. (s) & Reduction (s) \\
    \midrule
            ieee118 & 10  & {\bf 0.06}   & 16           & {\bf   0.17}  & 0.01  & 16           & {\bf   0.26}  & 0.01  & 0.02   \\
            ieee300 & 10  & {\bf 0.12}   & 22           & {\bf   0.28}  & 0.06  & 22           & {\bf   0.42}  & 0.02  & 0.03   \\
    ACTIVSg500      & 20  & {\bf 0.51}   & 24           & {\bf   0.29}  & 0.06  & 24           & {\bf   0.45}  & 0.02  & 0.03   \\
      1354pegase    & 22  & {\bf 1.17}   & 40           & {\bf   0.85}  & 0.35  & 40           & {\bf   1.16}  & 0.09  & 0.29   \\
    ACTIVSg2000     & 18  & {\bf 1.62}   & 43           & {\bf 1.95}    & 1.42  & 43           & {\bf 2.40}    & 0.24  & 0.65   \\
      2869pegase    & 22  & {\bf 2.11}   & 50           & {\bf   2.03}  & 1.18  & 50           & {\bf   2.69}  & 0.20  & 0.97   \\
      9241pegase    & 102 & {\bf  31.7}  & 69           & {\bf  10.65}  & 6.14  & 69           & {\bf  23.72}  & 1.17  & 16.23  \\
    ACTIVSg10k      & 130 & {\bf  39.3}  & 76$^\dagger$ & {\bf   7.9}   & 5.66  & 88$^\dagger$ & {\bf 21.9}    & 1.5   & 14.7   \\
    13659pegase     & 120 & {\bf  116.0} & 346          & {\bf   98.31} & 67.13 & 145          & {\bf 242.69}  & 19.23 & 202.89 \\
    ACTIVSg25k      & 47  & {\bf 36.1}   & 86           & {\bf  24.70}  & 16.86 & 86           & {\bf 84.96}   & 4.27  & 68.11  \\
    ACTIVSg70k      & 101 & {\bf 242.0}  & 90$^\dagger$ & {\bf  89.8}   & 65.7  & 85$^\dagger$ & {\bf 658.2}   & 21.5  & 606.5  \\
    \midrule
    9591goc         & 37  & {\bf 22.5}   & 43           & {\bf 11.66}   & 10.38 & 43           & {\bf   7.69}  & 2.13  & 1.61   \\
    10480goc        & 40  & {\bf 25.7}   & 50           & {\bf 13.98}   & 11.95 & 50           & {\bf   11.46} & 3.93  & 3.34   \\
    19402goc        & 45  & {\bf 66.5}   & 47           & {\bf 30.75}   & 26.83 & 47           & {\bf   19.52} & 4.86  & 7.24   \\
    \bottomrule
\end{tabular}
  }
  \caption{Benchmarking {\tt LinRed} with \emph{Full-space IPM}}
  \label{tab:results:benchmark}
\end{table}

We make the following observations.
(i) \emph{Hybrid full-space IPM} is slightly faster than Knitro, but only
because it evaluates the derivatives on the GPU.
(ii) {\tt LinRed} is able to solve all
the instances, including {\tt ACTIVSg70k}.
(iii) As expected, we get the same number of iterations
between \emph{Full-space IPM} and {\tt LinRed} except
on {\tt 13659pegase}. This case is indeed ill-conditioned, and the
convergence is sensitive to the linear solver employed (even MA27 and MA57 give
different results on this case).
(iv) Despite the good performance of the Cholesky factorization,
{\tt LinRed} is negatively impacted
by the scalability of the reduction algorithm: the larger the relative number of controls with respect to the total number of variables,
the less competitive is the reduction. On the largest instances,
{\tt LinRed} beats \emph{Full-space IPM} only
on the three {\tt goc} instances, which have fewer degrees of freedom
(see Table~\ref{tab:test_instances}).

\subsection{Linearize-then-reduce or reduce-then-linearize?}
\label{sec:results:redlin}
Now we benchmark {\tt LinRed} with its feasible
variant, {\tt RedLin}.
In  contrast to {\tt LinRed}, {\tt RedLin} follows a feasible path
w.r.t. the state equations~\eqref{eq:stateequation}: the algorithm can be stopped
at any time and return a feasible point, an advantageous feature if the resolution is time-constrained.
{\tt LinRed} uses the same setting as before,
and {\tt RedLin} is also using the MadNLP solver,
using the reduced derivatives defined in \eqref{eq:reduced:derivatives}.
{\tt RedLin} solves the state equations~\eqref{eq:stateequation}
(the power flow balance equations) at each iteration with a Newton--Raphson
algorithm, with a tolerance of $10^{-10}$.
The algorithm has two drawbacks: (i) this approach requires evaluating
the reduced Jacobian $\hat{A}_u$, with size $m \times n_u$, and
(ii) the default scaling computed by MadNLP
depends on $G_x^{-1}$, rendering the scaling inappropriate
if $G_x$ has a poor conditioning.
To ensure that the comparison is fair, we have modified our implementation
so that {\tt RedLin} uses the same scaling
as {\tt LinRed}.
The results are presented in Figure~\ref{fig:results:benchmark_reduced}
(the time spent in the reduction is omitted since it is the same as
in {\tt LinRed}).
\begin{figure}[!ht]
  \begin{subfigure}{.5\textwidth}
  \includegraphics[width=.9\textwidth]{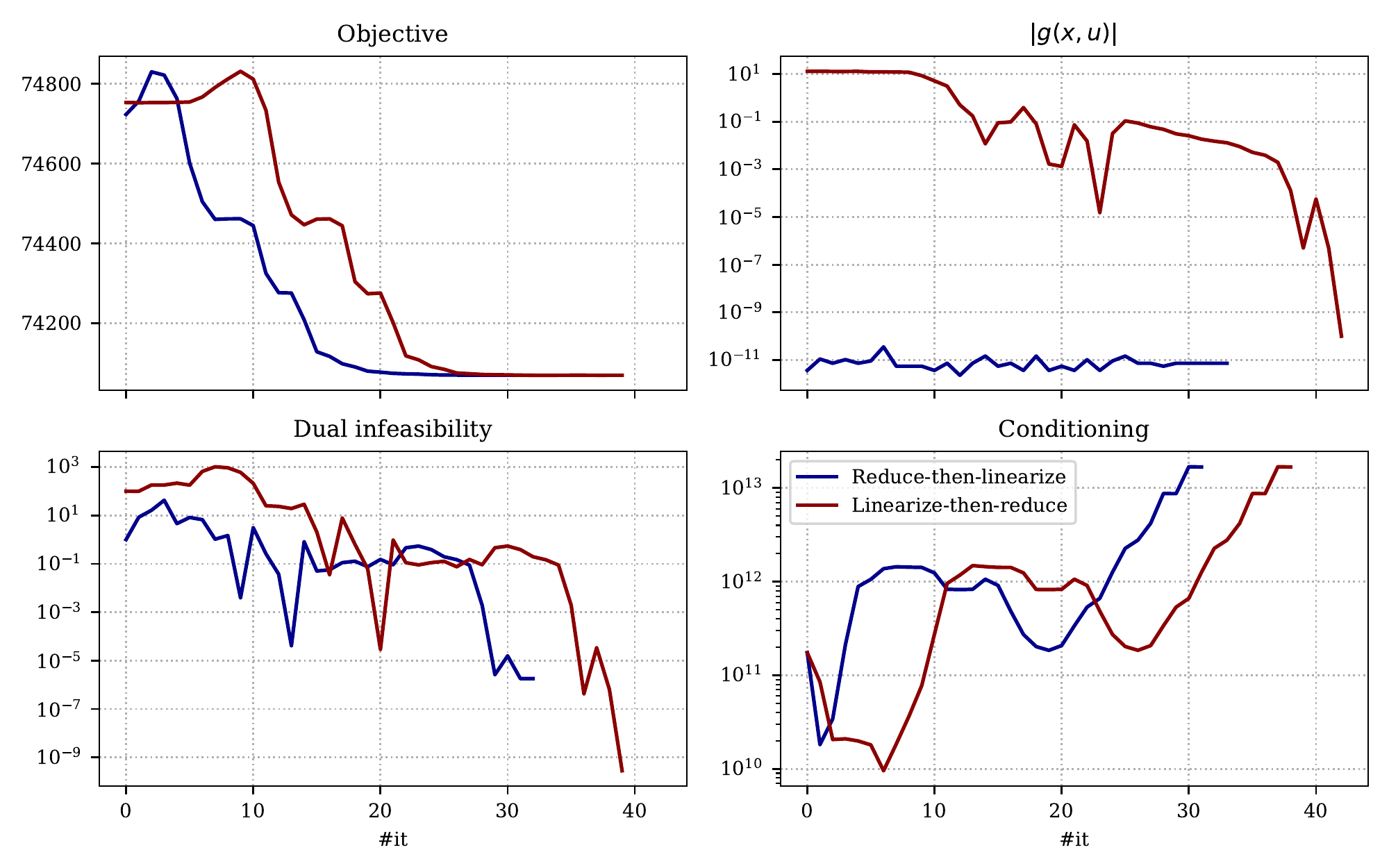}
  \caption{Convergence of {\tt RedLin} (blue)
    and {\tt LinRed} (red) on {\tt 1354pegase}.
  }
  \label{fig:results:convergence}
  \end{subfigure}
  \begin{subfigure}{.5\textwidth}
  \resizebox{.8\textwidth}{!}{
  \begin{tabular}{l|rrrr}
  \toprule
                  Case &   \#it &  Time (s) &  Chol. (s) &    PF (s) \\
  \midrule
             case118 & 16  & {\bf 0.40  } & 0.01  & 0.13  \\
             case300 & 24  & {\bf 0.55  } & 0.02  & 0.20  \\
     ACTIVSg500      & 22  & {\bf 0.58  } & 0.01  & 0.23  \\
      1354pegase     & 32  & {\bf 1.34  } & 0.05  & 0.35  \\
    ACTIVSg2000      & 33  & {\bf 2.25  } & 0.07  & 0.65  \\
      2869pegase     & 34  & {\bf 2.62  } & 0.08  & 0.67  \\
      9241pegase     & 48  & {\bf 28.30 } & 0.78  & 3.05  \\
     ACTIVSg10k      & 140 & {\bf 90.53 } & 3.60  & 5.51  \\
     13659pegase     & 167 & {\bf 944.04} & 33.19 & 15.78 \\
     ACTIVSg25k      & 52  & {\bf 156.17} & 1.79  & 5.36  \\
     ACTIVSg70k      & -   & -            & -     & -     \\
  \bottomrule
  \end{tabular}
  }
  \caption{Performance of {\tt RedLin}.}
  \end{subfigure}
  \caption{Benchmarking {\tt LinRed} with {\tt RedLin}.}
  \label{fig:results:benchmark_reduced}
\end{figure}

Comparing  with Table~\ref{tab:results:benchmark},
we make the following observations.
(i) {\tt RedLin} is able to solve instances
with up to 25,000 buses, which, to the best of our knowledge,
is a net improvement compared with previous attempts to solve
the OPF in the reduced-space~\cite{kardos2020reduced,pacaud2021feasible}.
(ii) On the largest instances, {\tt RedLin}
is penalized compared with {\tt LinRed}, since it has
to deal with the reduced Jacobian $\hat{A}_u$. For that reason,
{\tt RedLin} breaks on {\tt ACTIVSg70k}, since we
are running out of memory.
(iii) From {\tt case118} to {\tt 9241pegase},
{\tt RedLin} converges in fewer iterations
than does {\tt LinRed}.
(iv) The time to solve the power flow (column PF) is a fraction
of the time spent in the reduction algorithm.
(v) On {\tt ACTIVSg10k}, {\tt RedLin } does not
run into the feasibility restoration we encountered in {\tt LinRed}:
on this difficult instance, the convergence of {\tt RedLin} is smoother,
even if it requires more iterations.
(vi) the Newton--Raphson is not guaranteed to converge, but empirically we find this
is not an issue when a second-order method is employed.

In Figure~\ref{fig:results:convergence}
we compare the convergence of the two algorithms on {\tt 1354pegase}.
We observe that {\tt RedLin} converges faster than {\tt LinRed},
the latter satisfying the power flow equations only at the final iterations.
Interestingly, the conditioning of the KKT matrix $K_{cond}$ does not blow up
at the final iterations and remains below $10^{13}$, a reasonable value for an
interior-point algorithm.

\section{Conclusion}
This paper has presented an efficient implementation of the IPM
on GPU architectures, based on a Schur reduction of the underlying KKT system.
We have derived two practical algorithms, \emph{linearize-then-reduce} and
\emph{reduce-then-linearize}, adapted their workflows to be efficient when
the vast majority of their computation is run on the GPU, and detailed their respective performance
on different large-scale OPF instances.
We also discussed the benefits of a hybrid full space IPM solver --- computing the
derivatives on the GPU and the linear algebra on the CPU --- and demonstrated
that this approach generally outperforms Knitro, running exclusively on the CPU.
The relative performance of the reduced-space algorithms is highly dependent on
the ratio of controls with respect to the total number of variables.
Their performance lags behind both
Knitro and the hybrid full space solver when the problem has many
control variables (as it is the case on the MATPOWER benchmark)
but is significantly ahead -- up to a factor of 3 -- when
the problem has a relatively lower number of control variables.
Moreover, the reduce-then-linearize algorithm has the added
benefit of producing a feasible solution to the power flow equations at
\emph{any} iteration, which makes it a great candidate for real time applications.
To improve the performance of reduction algorithms, we believe the most
important item is to alleviate the dependence on the
number of control variables. We plan to explore a way to accelerate the
reduction by exploiting the exponentially decaying structure of the reduced
Hessian~\cite{shin2021exponential}.

\section*{Acknowledgments}
We thank
Juraj Kardo\v{s} for reviving the interest in reduced-space methods and
Kasia \'{S}wirydowicz for pointing us to the {\tt cusolverRF} library.
This research was supported by the Exascale Computing Project (17-SC-20-SC),
a joint project of the U.S. Department of Energy’s Office of Science and
National Nuclear Security Administration, responsible for delivering a
capable exascale ecosystem, including software, applications, and hardware
technology, to support the nation’s exascale computing imperative.
The material is based upon work supported in part by the U.S. Department of Energy, Office of Science, under contract  DE-AC02-06CH11357.


\vfill
\begin{flushright}
{\footnotesize
  \framebox{\parbox{\textwidth}{
The submitted manuscript has been created by UChicago Argonne, LLC,
Operator of Argonne National Laboratory (``Argonne"). Argonne, a
U.S. Department of Energy Office of Science laboratory, is operated
under Contract No. DE-AC02-06CH11357. The U.S. Government retains for
itself, and others acting on its behalf, a paid-up nonexclusive,
irrevocable worldwide license in said article to reproduce, prepare
derivative works, distribute copies to the public, and perform
publicly and display publicly, by or on behalf of the Government.
The Department of
Energy will provide public access to these results of federally sponsored research in accordance
with the DOE Public Access Plan. http://energy.gov/downloads/doe-public-access-plan. }}
\normalsize
}
\end{flushright}

\end{document}